\documentclass[11pt]{amsart}

\usepackage{amsfonts,amssymb,amsmath,amsthm}
\usepackage{bbm}
\usepackage{hhline}
\usepackage{stmaryrd}
\usepackage{longtable}
\usepackage{fancyhdr}

\input{xy}
\xyoption{all}

\begin{document}

\setlength{\parindent}{0pt} \setlength{\parskip}{1ex}

\newtheorem{Le}{Lemma}[section]
\newtheorem{Pro}[Le]{Proposition}
\newtheorem{Th}[Le]{Theorem}
\newtheorem{Co}[Le]{Corollary}
\newtheorem{Conj}[Le]{Conjecture}

\newtheorem{Thm}{Theorem}

\theoremstyle{definition}
\newtheorem{Def}[Le]{Definition}
\newtheorem{DaL}[Le]{Definition and Lemma}
\newtheorem{Conv}[Le]{Convention}
\newtheorem{Assumpt}[Le]{Assumption}

\theoremstyle{remark}
\newtheorem{Ex}[Le]{Example}
\newtheorem{Rem}[Le]{Remark}

\author{Frank Reidegeld}
\title{Spaces admitting homogeneous $G_2$-structures}
\date{}
\address{E-Mail: frank.reidegeld\char"40 math.uni-hamburg.de}

\begin{abstract}
We classify all seven-dimensional spaces which admit
a ho\-mogeneous cosymplectic $G_2$-structure. The motivation for
this classification is that each of these spaces is a possible
principal orbit of a parallel $\text{Spin}(7)$-manifold of
cohomogeneity one.
\end{abstract}

\maketitle

\section{Introduction}

The aim of this article is to classify all spaces which admit a
homogeneous cosymplectic $G_2$-structure. Moreover, we not only
classify the spaces themselves, but also the transitive group
actions which preserve at least one cosymplectic $G_2$-structure.

In the literature, many of such spaces are known. Friedrich, Kath,
Moroianu, and Semmelmann \cite{Fried} classify all simply
connected, compact spaces which admit a homogeneous nearly parallel
$G_2$-structure. The product of a space with a
homogeneous $SU(3)$-structure and a circle carries a canonical
homogeneous $G_2$-structure. The spaces from the article of Cleyton
and Swann \cite{Cley} which admit a homogeneous $SU(3)$-structure
should therefore be mentioned in this context, too.

One reason for our interest in this kind of spaces is that
any principal orbit of a parallel $\text{Spin}(7)$-manifold of
cohomogeneity one carries a homogeneous cosymplectic 
$G_2$-structure. Conversely, any homogeneous cosymplectic 
$G_2$-structure can be extended to a parallel
$\text{Spin}(7)$-manifold of cohomogeneity one. A discussion of
these facts can be found in Hitchin \cite{Hitchin}. The aim of
this article is to prove the following theorem:

\begin{Thm} \label{MainThm}
\begin{enumerate}
    \item \label{RedTh}
    Let $G/H$ be a seven-dimensional, compact, connected, 
    $G$-homogeneous space which admits a $G$-invariant $G_2$-structure.
    We assume that $G/H$ is a product of a circle and another 
    homogeneous space and that $G$ acts almost effectively on $G/H$. 
    Furthermore, we assume that $G$ and $H$ are both connected. In this
    situation, $G$, $H$, and $G/H$ are up to a covering one of the 
    spaces from the table below:

    \newpage

    \begin{center}
    \begin{longtable}{|l|l|l|}
    \hline $G$ & $H$ & $G/H$ \\

    \hline \hline\endhead

    $U(1)^7$ & $\{e\}$ & $T^7$ \\

    \hline

    $SU(2)\times U(1)^4$ & $\{e\}$ & $S^3\times T^4$ \\

    \hline

    $SU(2)^2\times U(1)$ & $\{e\}$ & $S^3\times S^3\times S^1$ \\

    \hline

    $SU(2)^2\times U(1)^2$ & $U(1)$ & $S^3\times S^3\times S^1$ \\

    \hline

    $SU(2)^2\times U(1)^2$ & $U(1)$ & $SU(2)^2/U(1)\times T^2$ \\

    \hline

    $SU(2)^3\times U(1)$ & $SU(2)$ & $S^3\times S^3\times S^1$ \\

    \hline

    $SU(3)\times U(1)^2$ & $SU(2)$ & $S^5\times T^2$ \\

    \hline

    $SU(3)\times U(1)$ & $U(1)^2$ & $SU(3)/U(1)^2\times S^1$ \\

    \hline

    $Sp(2)\times U(1)$ & $Sp(1)\times U(1)$ & $\mathbb{CP}^3\times
    S^1$ \\

    \hline

    $G_2\times U(1)$ & $SU(3)$ & $S^6\times S^1$ \\

    \hline
    \end{longtable}
    \end{center}
    Conversely, any of the above spaces admits a $G$-invariant
    $G_2$-structure.
   
    \item \label{IrredTh}
    Let $G$, $H$, and $G/H$ satisfy the same conditions as in 
    (\ref{RedTh}) with the single exception that $G/H$ is 
    \underline{not} a product of a circle and another homogeneous
    space. In this situation, $G$, $H$, and $G/H$ are up to a 
    covering one of the spaces from the table below:

    \begin{center}
    \begin{longtable}{|l|l|l|}
    \hline

    $G$ & $H$ & $G/H$ \\

    \hline \hline

    $SU(3)$ & $U(1)$ & $N^{k,l}\quad\text{with}\:\:k,l\in\mathbb{Z}$
    \\

    \hline

    $SO(5)$ & $SO(3)$ & $V^{5,2}$ \\

    \hline

    $Sp(2)$ & $Sp(1)$ & $S^7$ \\

    \hline

    $SO(5)$ & $SO(3)$ & $B^7$ \\

    \hline

    $SU(2)^3$ & $U(1)^2$ & $Q^{1,1,1}$ \\

    \hline

    $SU(3)\times U(1)$ & $U(1)^2$ &
    $N^{k,l}\quad\text{with}\:\:k,l\in\mathbb{Z}$ \\

    \hline

    $SU(3)\times SU(2)$ & $SU(2)\times U(1)$ & $M^{1,1,0}$ \\

    \hline

    $SU(3)\times SU(2)$ & $SU(2)\times U(1)$ & $N^{1,1}$ \\

    \hline

    $Sp(2)\times U(1)$ & $Sp(1)\times U(1)$ & $S^7$ \\

    \hline

    $Sp(2)\times Sp(1)$ & $Sp(1)\times Sp(1)$ & $S^7$ \\

    \hline

    $SU(4)$ & $SU(3)$ & $S^7$ \\

    \hline

    $\text{Spin}(7)$ & $G_2$ & $S^7$ \\

    \hline
    \end{longtable}
    \end{center}
    Conversely, any of the above spaces admits a $G$-invariant
    $G_2$-struc\-ture.
    
    \item Any of the spaces $G/H$ from (\ref{RedTh}) or (\ref{IrredTh})
    even admits a $G$-invariant cosymplectic $G_2$-structure.
\end{enumerate}
\end{Thm}

In table (\ref{IrredTh}), $N^{k,l}$ denotes an Aloff-Wallach space,
$V^{5,2}$ denotes the Stiefel manifold of all orthonormal pairs in
$\mathbb{R}^5$, and $B^7$ is the seven-dimensional Berger space. In
the fourth, fifth, and sixth row of table (\ref{RedTh}) and in the fifth 
and seventh row of table (\ref{IrredTh}), the embedding of $H$ into $G$ has
to be special in order to make $G/H$ a space which admits a $G$-invariant
$G_2$-structure. The details of those embeddings are described in Section 
\ref{RedSec} and \ref{IrredSec}. In the other cases, the information in the
above tables is sufficient to determine the embedding of $H$ into $G$.

From the theorem it follows that either $G/H$ is a product of a circle and
a space which admits a homogeneous $SU(3)$-structure or that it cannot be
decomposed into factors of lower dimension. We remark that we not only prove
the existence of a homogeneous cosymplectic $G_2$-structure on each of the 
spaces but also the existence of cosymplectic $G_2$-structures which are 
invariant under any of the transitive group actions. The space $(SU(2)
\times SU(2))/U(1) \times T^2$ admits a homogeneous $G_2$-structure but
seems not to be mentioned in the literature before.

The proof of Theorem \ref{MainThm} consists of three steps: After two
introductory sections, we classify all connected Lie subgroups of $G_2$.
This is necessary, since in the situation of the theorem $H$ can be embedded into
$G_2$. In Section \ref{RedSec} and \ref{IrredSec}, we determine all $G/H$ which
admit a $G$-invariant, but not necessarily cosymplectic $G_2$-structure.
Finally, we have to prove the existence of a $G$-invariant cosymplectic
$G_2$-structure on all of the spaces which we have found. This will be done
in Section \ref{CosympSec}.

\section{The group $G_2$}
\label{2ndSec}

Before we classify the connected subgroups of $G_2$, we collect some facts
on this group. For a more comprehensive introduction into this issue, see 
Baez \cite{Baez} or Bryant \cite{Br}.

The group $G_2$ can be defined with help of the octonions: We
recall that a \textit{normed division algebra} is a pair
$(A,\langle\cdot,\cdot\rangle)$ of a real, not necessarily
associative algebra with a unit element and a scalar product which
satisfies $\langle x\cdot y,x\cdot y\rangle = \langle x,x\rangle\:
\langle y,y\rangle$ for all $x,y\in A$. There exists up to
isomorphisms exactly one eight-dimensional normed division
algebra, namely the \textit{octonions} $\mathbb{O}$.

The quaternions $\mathbb{H}$ are a subalgebra of $\mathbb{O}$. We
fix an octonion $\epsilon$ in the orthogonal complement of
$\mathbb{H}$ such that $\|\epsilon\|=1$. We call
$(x_0,\ldots,x_7):=(1,i,j,k,\epsilon,i\epsilon,
j\epsilon,k\epsilon)$ the \textit{standard basis of $\mathbb{O}$}.
Let $\text{Im}(\mathbb{O}):=\text{span}(1)^\perp$ be the
\textit{imaginary space of $\mathbb{O}$}. The map

\begin{equation}
\begin{split}
\omega &
:\text{Im}(\mathbb{O})\times\text{Im}(\mathbb{O})\times\text{Im}(
\mathbb{O})\rightarrow\mathbb{R}\\ \omega(x,y,z) & :=\langle
x\cdot y,z\rangle\\
\end{split}
\end{equation}

is a three-form. From now on, we denote $dx^{i_1}\wedge \ldots
\wedge dx^{i_k}$ shortly by $dx^{i_1\ldots i_k}$. With this
notation, we have:

\begin{equation}
\label{omega}
\omega=dx^{123}+dx^{145}-dx^{167}+dx^{246}+dx^{257}+dx^{347}-dx^{356}\:.
\end{equation}

\begin{Rem}
The multiplication table of $\mathbb{O}$ is uniquely determined by
the coefficients of $\omega$. Let $\epsilon'$ be an octonion with
the same properties as $\epsilon$. Since there exists an
automorphism of $\mathbb{O}$ which is the identity on $\mathbb{H}$
and maps $\epsilon$ to $\epsilon'$, $\omega$ is independent of the
choice of $\epsilon$.
\end{Rem}

We are now able to define the Lie group $G_2$:

\begin{DaL}
\begin{enumerate}
    \item Any automorphism $\varphi$ of $\mathbb{O}$ satisfies
    $\varphi(\text{Im}(\mathbb{O})) \subseteq$ $\text{Im}(\mathbb{O})$ and thus
    can be identified with a map from $\text{Im}(\mathbb{O})$ onto itself. $G_2$
    is defined as the stabilizer group of $\omega$ or equivalently as the
    automorphism group of $\mathbb{O}$.
    \item The Lie algebra of $G_2$ we denote by $\mathfrak{g}_2$.
    \item The seven-dimensional representation which is induced by the action of
    $G_2$ on $\text{Im}(\mathbb{O})$ by automorphisms we call the
    \textit{standard representation of $G_2$}.
\end{enumerate}
\end{DaL}

A proof of the fact that the stabilizer of $\omega$ is the same as
the automorphism group of $\mathbb{O}$ can be found in Bryant
\cite{Br}. Later on, we work with the Hodge dual
$\ast\omega\in\bigwedge^4 \text{Im}( \mathbb{O})^\ast$ of $\omega$
which is taken with respect to $\langle\cdot, \cdot\rangle$ and
the orientation which makes $(x_1,\ldots,x_7)$ positive:

\begin{equation}
\label{astomega}
\ast\omega=-dx^{1247}+dx^{1256}+dx^{1346}+dx^{1357}-dx^{2345}+dx^{2367}
+dx^{4567}\:.
\end{equation}

Finally, we fix a Cartan subalgebra $\mathfrak{t}$ of
$\mathfrak{g}_2$, which we will need for our explicit
calculations. With respect to the standard basis of
$\text{Im}(\mathbb{O})$, $\mathfrak{t}$ is the following set of
matrices:

\begin{equation}
\label{csag2} \mathfrak{t}:= \left\{\left(\,\begin{array}{ccccccc}
\hhline{-~~~~~~} \multicolumn{1}{|c|}{0}\\ \hhline{---~~~~} &
\multicolumn{1}{|c}{0} & \multicolumn{1}{c|}{\lambda_1} &&&&\\ &
\multicolumn{1}{|c}{-\lambda_1} & \multicolumn{1}{c|}{0} &&&&\\
\hhline{~----~~} &&& \multicolumn{1}{|c}{0} &
\multicolumn{1}{c|}{\lambda_2} &&\\ &&&
\multicolumn{1}{|c}{-\lambda_2} & \multicolumn{1}{c|}{0} &&\\
\hhline{~~~----} &&&&& \multicolumn{1}{|c}{0} &
\multicolumn{1}{c|}{\lambda_1+\lambda_2}\\ &&&&&
\multicolumn{1}{|c}{-\lambda_1-\lambda_2} &
\multicolumn{1}{c|}{0}\\ \hhline{~~~~~--}
\end{array}\,\right)\:\left|\begin{array}{c} \\ \\ \\ \\ \\ \\ \\
\end{array}\right.\!\!\!\!\!
\lambda_1,\lambda_2\in\mathbb{R}\right\}
\end{equation}

\section{Some remarks on $G_2$-structures}

In this section, we introduce the different types of
$G_2$-structures which we consider in this article. We refer the
reader to Bryant \cite{Br} or the books of Joyce \cite{Joy} 
and Salamon \cite{Sal} for further facts on these structures. A 
$G_2$-structure can be defined as a three-form which is at each point
stabilized by $G_2$:

\begin{Def}
Let $M$ be a seven-dimensional manifold and $\omega$ be a
three-form on $M$ with the following property: For any $p\in M$
there exists a neighborhood $U$ of $p$ and vector fields
$X_1,\ldots,X_7$ on $U$ such that

\begin{equation}
\omega_q(X_i,X_j,X_k)=\omega(x_i,x_j,x_k)\quad\forall q\in U,\:
i,j,k\in\{1,\ldots,7\}\:.
\end{equation}

The $\omega$ on the right hand side of the above formula is the
three-form (\ref{omega}) and $x_i$, $x_j$, and $x_k$ are elements
of the standard basis of $\mathbb{O}$. In this situation,
$\omega$ is called a \textit{$G_2$-structure on $M$} and the pair
$(M,\omega)$ is called a \textit{$G_2$-manifold}.
\end{Def}

On any $G_2$-manifold $(M,\omega)$ there exist a metric $g$ and a
volume form $\text{vol}$ which are defined by:

\begin{equation}
g(X,Y)\:\text{vol}:=-\frac{1}{6}\:(X\:\lfloor\:\omega) \wedge
(Y\:\lfloor\:\omega)\wedge\omega\:.
\end{equation}

We call $g$ the \textit{associated metric} and $\text{vol}$ the
\textit{associated volume form}. $g$ and $\text{vol}$ induce a
Hodge star operator $\ast:\bigwedge^{\ast} T^\ast M \rightarrow
\bigwedge^\ast T^\ast M$ and we therefore have a four-form
$\ast\omega$ on $M$, which is invariant under the stabilizer $G_2$
of $\omega$. On the flat $G_2$-manifold $(\mathbb{R}^7, \omega)$
this four-form coincides with (\ref{astomega}). For our
considerations, we need the following types of $G_2$-structures:

\begin{Def}
A $G_2$-manifold $(M,\omega)$ is called

\begin{enumerate}
    \item \textit{parallel} if $d\omega=0$ and $d\ast\omega=0$,
    \item \textit{nearly parallel} if there exists a $\lambda\in\mathbb{R}
    \setminus\{0\}$ such that $d\omega=\lambda\ast\omega$ and thus $d\ast\omega=0$,
    \item \textit{cosymplectic} if $d\ast\omega=0$.
\end{enumerate}
\end{Def}

Further information on the different types of $G_2$-structures can
be found in the article by Fern\'andez and Gray \cite{Fer}. We will
deal first of all with homogeneous $G_2$-manifolds:

\begin{Def}
A $G_2$-manifold $(M,\omega)$ is called \textit{($G$-)homogeneous} if
there exists a transitive smooth action by a Lie group $G$ which
leaves $\omega$ invariant.
\end{Def}

In the above situation, $M$ is $G$-equivariantly diffeomorphic to a
quotient $G/H$. $G$ can be chosen in such a way that it acts effectively on
$G/H$ and preserves $\omega$. $H$ acts on the tangent space of
$G/H$ by its isotropy representation. Since $G_2$ acts on the
tangent space as the stabilizer of $\omega$ and $\omega$ is
$G$-invariant, we have proven the following lemma:

\begin{Le} \label{G2Hom1}
Let $G/H$ be a seven-dimensional $G$-homogeneous space which admits a
$G$-invariant $G_2$-structure. We assume that $G$ acts effectively
on $G/H$. In this situation, there exists a vector space
isomorphism $\varphi:T_p G/H\rightarrow \mathbb{R}^7$ such that
$\varphi H \varphi^{-1}\subseteq G_2$, where $H$ is identified
with its isotropy representation and $G_2$ with its
seven-dimensional irreducible representation.
\end{Le}

The converse of the above lemma is also true:

\begin{Le} \label{G2Hom2}
Let $G/H$ be a seven-dimensional $G$-homogeneous space such
that $G$ acts effectively and there exists a vector space isomorphism
$\varphi:T_p G/H\rightarrow \mathbb{R}^7$ with $\varphi H
\varphi^{-1}\subseteq G_2$. In this situation, there exists a
$G$-invariant $G_2$-structure on $G/H$.
\end{Le}

$\mathbf{Proof:}$ The action of $G$ on the tangent bundle
determines a $G$-invariant $H$-structure on $G/H$. Its extension
to a principal bundle with structure group $G_2$ is a
$G$-invariant $G_2$-structure.

\begin{flushright}
$\Box$
\end{flushright}

\section{Subgroups of $G_2$}
\label{G2Sec}

In this section, we classify all connected subgroups of $G_2$.
First, we describe all of these subgroups explicitly. After that,
we prove that the list which we have found is complete.

On page \pageref{csag2}, we have described a Cartan subalgebra
$\mathfrak{t}$ of $\mathfrak{g}_2$. $\mathfrak{t}$ is generated by
the two elements which satisfy $(\lambda_1,\lambda_2)=(1,0)$ and
$(\lambda_1,\lambda_2)=(0,1)$. $\text{Im}(\mathbb{O})$ splits with
respect to the action of $\mathfrak{t}$ into
$\mathbbm{V}_{1,0}^\mathbb{C} \oplus
\mathbbm{V}_{0,1}^\mathbb{C}\oplus \mathbbm{V}_{1,1}^\mathbb{C}
\oplus \mathbbm{V}_{0,0}^\mathbb{R}$. The subscripts denote the
weights with which the two generators of $\mathfrak{t}$ act and
the superscript indicates if the submodule is complex or real.
Since any abelian subalgebra of $\mathfrak{g}_2$ is conjugate to a
subalgebra of $\mathfrak{t}$, we have finished the abelian case.

Next, we describe the subgroups of $G_2$ whose Lie algebra has an
ideal of type $\mathfrak{su}(2)$. In an article by Cacciatori et
al. \cite{Cac}, the authors introduce the following Lie group
homomorphism:

\begin{equation}
\begin{split}
\varphi & : Sp(1)\times Sp(1) \rightarrow G_2\\
\varphi(h,k)(x+y\epsilon) & :=hxh^{-1} + (kyh^{-1})\epsilon\:,\\
\end{split}
\end{equation}

where $x,y\in\mathbb{H}$ and $Sp(1)$, which is isomorphic to $SU(2)$,
is identified with the unit quaternions. The kernel of $\varphi$ is
$\{(1,1),(-1,-1)\}$ and its image thus is isomorphic to $SO(4)$.
The first factor of $Sp(1)\times Sp(1)$ acts irreducibly on
$\text{Im}(\mathbb{H})$ and $\mathbb{H}\epsilon$ and the second
factor acts irreducibly on $\mathbb{H}\epsilon$ and trivially on
its orthogonal complement. The splitting of
$\text{Im}(\mathbb{O})$ into irreducible $2\mathfrak{su}(2)$-modules 
therefore is $\mathbbm{V}_{2,0}^{\mathbb{R}} \oplus 
\mathbbm{V}_{1,1}^{\mathbb{C}}$. Analogously to above, the subscripts
of the modules denote the weights of the $2\mathfrak{su}(2)$-action
with respect to the first and second summand. By a straightforward
calculation, we can prove that the group $Sp(1)$ which is
diagonally embedded into $Sp(1)\times Sp(1)$ acts irreducibly on
$\text{Im}(\mathbb{H})$ and $\text{Im}(\mathbb{H})\epsilon$ and
trivially on $\text{span}(\epsilon)$. $\text{Im}(\mathbb{O})$ thus
splits into $2\mathbbm{V}_2^{\mathbb{R}} \oplus
\mathbbm{V}_0^{\mathbb{R}}$ with respect to that subgroup.

In his article "Semisimple subalgebras of semisimple Lie algebras"
\cite{Dynkin}, Dynkin proves the existence of another subalgebra
of $\mathfrak{g}_2$ which is isomorphic to $\mathfrak{su}(2)$ and
acts irreducibly on $\text{Im}(\mathbb{O})$ with weight $6$. Since
we do not need an explicit description of that subalgebra, we
simply state its existence.

According to the non zero weights of their action on
$\text{Im}(\mathbb{O})$, we denote the four subalgebras of
$\mathfrak{g}_2$ which are isomorphic to $\mathfrak{su}(2)$ by
$\mathfrak{su}(2)_{1,2}$, $\mathfrak{su}(2)_1$,
$\mathfrak{su}(2)_{2,2}$, and $\mathfrak{su}(2)_6$.

By a short calculation, we see that the element of $\mathfrak{t}$
with $\lambda_1=2$ and $\lambda_2=-1$ commutes with
$\mathfrak{su}(2)_1$ and the element with $\lambda_1=0$ and
$\lambda_2=1$ commutes with $\mathfrak{su}(2)_{1,2}$.
$\mathfrak{g}_2$ therefore contains a subalgebra of type
$\mathfrak{su}(2)_1 \oplus \mathfrak{u}(1)$ and a subalgebra of
type $\mathfrak{su}(2)_{1,2} \oplus \mathfrak{u}(1)$. Both of them
are a direct sum of an ideal of $\mathfrak{su}(2) \oplus
\mathfrak{su}(2) \subseteq \mathfrak{g}_2$ and a one-dimensional
subalgebra of the other ideal.

The group of all automorphisms of $\mathbb{O}$ which fix $i$ is a
compact, connected, eight-dimensional Lie group. Its action on
$\mathbb{C}$ is trivial and it acts irreducibly on the orthogonal
complement of $\mathbb{C}\subseteq \mathbb{O}$. These conditions
force the group to be isomorphic to $SU(3)$.

Our next step is to prove that there are up to conjugation by an
element of $G_2$ no further connected subgroups. $\mathfrak{g}_2$
is a Lie algebra of rank $2$ and two maximal tori of a Lie group
are always conjugate to each other. Therefore, further
connected, abelian subgroups of $G_2$ cannot exist.

According to Dynkin \cite{Dynkin}, all semisimple subalgebras of
$\mathfrak{g}_2$ are isomorphic to $\mathfrak{su}(2)$,
$2\mathfrak{su}(2)$, $\mathfrak{su}(3)$, or $\mathfrak{g}_2$. Any
of these algebras acts by the restriction of the adjoint
representation on $\mathfrak{g}_2$. The weights of this action are
computed in \cite{Dynkin}, too. The list Dynkin obtains is the
same as our list of semisimple Lie subalgebras. Moreover, Dynkin
\cite{Dynkin} proves that his list is complete up to conjugation
by elements of $G_2$.

It remains to prove that there are no further subalgebras of type
$\mathfrak{su}(2)\oplus \mathfrak{u}(1)$. Let $x$ be a generator
of the center. Since $\mathfrak{su}(2)$ commutes with
$\mathfrak{u}(1)$, the action of $x$ on $\text{Im}(\mathbb{O})$
has to be $\mathfrak{su}(2)$-equivariant. With help of the real
version of Schur's Lemma, we are able to classify all
$\mathfrak{su}(2)$-equivariant endomorphisms of 
$\text{Im}(\mathbb{O})$ for any embedding of
$\mathfrak{su}(2)$ into $\mathfrak{g}_2$. Since the action of $x$
on $\text{Im}(\mathbb{O})$ has to be a restricted automorphism of
$\mathbb{O}$, we can reduce the number of those endomorphisms even further.
After that, we see that all $x\in\mathfrak{g}_2$ which commute with 
$\mathfrak{su}(2)_1$ ($\mathfrak{su}(2)_{1,2}$) are conjugate to the two
matrices which we have already found. The conjugation is with respect to an
element of the Lie subgroup of $G_2$ which is associated to
$\mathfrak{su}(2)_{1,2}$ ($\mathfrak{su}(2)_1$). By the same
method, we are able to prove that no subalgebras of type
$\mathfrak{su}(2)_{2,2}\oplus \mathfrak{u}(1)$ or
$\mathfrak{su}(2)_6 \oplus \mathfrak{u}(1)$ do exist. We finally
have proven the following theorem:

\begin{Th} \label{G2SubThm} Let $H$ be a connected Lie subgroup of $G_2$. We
denote the Lie algebra of $H$ by $\mathfrak{h}$. The irreducible
action of $G_2$ on $\text{Im}(\mathbb{O})$ induces an action of
$H$ on $\text{Im}(\mathbb{O})$. In this situation, $\mathfrak{h}$,
$H$, and the action of $H$ on $\text{Im}(\mathbb{O})$ are
contained in the table below. Moreover, any two connected Lie
subgroups of $G_2$ whose action on $\text{Im}(\mathbb{O})$ is
equivalent are conjugate not only by an element of $GL(7)$ but even
by an element of $G_2$.

\begin{center}
\begin{longtable}{|l|l|l|}
\hline

$\mathfrak{h}$ & $H$ & Splitting of $\text{Im}(\mathbb{O})$ into
irreducible summands\\

\hline \hline\endhead

$\{0\}$ & $\{e\}$ & \\

\hline

$\mathfrak{u}(1)$ & $U(1)$ & $\mathbbm{V}_a^{\mathbb{C}}\oplus
\mathbbm{V}_b^{\mathbb{C}} \oplus
\mathbbm{V}_{-a-b}^{\mathbb{C}}\oplus \mathbbm{V}_0^{\mathbb{R}}$
\\

\hline

$2\mathfrak{u}(1)$ & $U(1)^2$ & $\mathbbm{V}_{1,0}^\mathbb{C}
\oplus \mathbbm{V}_{0,1}^\mathbb{C}\oplus
\mathbbm{V}_{1,1}^\mathbb{C} \oplus \mathbbm{V}_{0,0}^\mathbb{R}$
\\

\hline

$\mathfrak{su}(2)$ & $SU(2)$ & $\mathbbm{V}_1^{\mathbb{C}}\oplus 3
\mathbbm{V}_0^{\mathbb{R}}$ \\

\hline

$\mathfrak{su}(2)$ & $SU(2)$ & $\mathbbm{V}_2^{\mathbb{R}}\oplus
\mathbbm{V}_1^{\mathbb{C}}$ \\

\hline

$\mathfrak{su}(2)$ & $SO(3)$ & $2\mathbbm{V}_2^{\mathbb{R}}\oplus
\mathbbm{V}_0^{\mathbb{R}}$ \\

\hline

$\mathfrak{su}(2)$ & $SO(3)$ & $\mathbbm{V}_6^{\mathbb{R}}$ \\

\hline

$\mathfrak{su}(2)\oplus\mathfrak{u}(1)$ & $U(2)$ &
$\mathbbm{V}_1^{\mathbb{C}}\oplus 3 \mathbbm{V}_0^{\mathbb{R}}$
w.r.t. $\mathfrak{su}(2)$ \\

\hline

$\mathfrak{su}(2)\oplus\mathfrak{u}(1)$ & $U(2)$ &
$\mathbbm{V}_2^{\mathbb{R}}\oplus \mathbbm{V}_1^{\mathbb{C}}$
w.r.t $\mathfrak{su}(2)$ \\

\hline

$2\mathfrak{su}(2)$ & $SO(4)$ & $\mathbbm{V}_{2,0}^{\mathbb{R}}
\oplus \mathbbm{V}_{1,1}^{\mathbb{C}}$\\

\hline

$\mathfrak{su}(3)$ & $SU(3)$ &
$\mathbbm{V}_{1,0}^{\mathbb{C}}\oplus
\mathbbm{V}_{0,0}^{\mathbb{R}}$
\\

\hline

$\mathfrak{g}_2$ & $G_2$ & $\mathbbm{V}_{1,0}^{\mathbb{R}}$ \\

\hline
\end{longtable}
\end{center}

The subscripts of the modules in the above table denote the
weights of the $H$-action and the superscript indicates if the
module is complex or real. Further details of the embeddings, in
particular of those of $U(2)$ and $SO(4)$ into $G_2$, we have
described on the preceding pages.
\end{Th}

\begin{Rem}
Most statements of Theorem \ref{G2SubThm} can also be proven by
elementary calculations which make use of the octonions. In order to keep our
presentation of this issue short, we often made use of the results
of Dynkin \cite{Dynkin}.
\end{Rem}

\section{The reducible case}
\label{RedSec}

We divide the spaces which admit a homogeneous $G_2$-structure
into two classes:

\begin{Def}
Let $G/H$ be a $G$-homogeneous space. We call $G/H$ \textit{$S^1$-reducible} if
it is $G$-equivariantly covered by a product of a circle and another homogeneous 
space. Otherwise, $G/H$ is called \textit{$S^1$-irreducible}.
\end{Def}

In this section, we classify all $S^1$-reducible spaces which admit a
homogeneous $G_2$-structure, and in the next section, we classify
the $S^1$-irreducible ones. We will see that none of the $S^1$-irreducible
spaces is covered by a product of lower-dimensional homogeneous spaces. 
The $S^1$-irreducible spaces which we will find are thus irreducible in 
the classical sense, too.

Throughout this article we denote the Lie algebra of $G$ by
$\mathfrak{g}$ and the Lie algebra of $H$ by $\mathfrak{h}$. In
order to simplify our considerations, we assume that $G/H$ is
compact and that $G$ is connected and acts almost effectively on
$G/H$, i.e. the subgroup of $G$ which acts as the identity map is
finite. Moreover, we classify the possible $G/H$ and $G$ only up
to coverings. Before we start our classification, we collect some
helpful facts:

\begin{enumerate}
    \item We have $\dim{\mathfrak{g}}=\dim{\mathfrak{h}} + 7$. This fact reduces
    the number of possible $\mathfrak{g}$ which we have to consider.

    \item Since $G/H$ is compact and $G$ is a subgroup of the isometry group of
    the metric on $G/H$, $G$ has to be compact, too. We thus can assume that
    $\mathfrak{g}$ is the direct sum of a semisimple and an abelian Lie algebra.

    \item Since the roots of a semisimple Lie algebra are paired, we
    have $\dim{\mathfrak{k}} \equiv \text{rank}\:{\mathfrak{k}} \:\:\: (mod\:\:\: 2)$ for
    any Lie algebra $\mathfrak{k}$ of a compact Lie group. It follows from 
    $\dim{\mathfrak{g}} = \dim{\mathfrak{h}} + 7$ that $\text{rank}\:{\mathfrak{g}}
    \not\equiv \text{rank}\:{\mathfrak{h}} \:\:\: (mod\:\:\: 2)$. $\mathfrak{h}$ can
    be considered as a subalgebra of $\mathfrak{g}_2$ and thus is trivial or of rank
    $1$ or $2$. 
    
    \begin{enumerate}
        \item If $\mathfrak{h}$ is of rank $1$, $\text{rank}\:{\mathfrak{g}}$ has
        to be even. The Cartan subalgebra of $\mathfrak{h}$ has to act on
        the tangent space in the same way as a one-dimensional subalgebra of 
        $\mathfrak{t}$ on $\text{Im}(\mathbb{O})$. The maximal trivial 
        $\mathfrak{h}$-submodule of the tangent space therefore is at most 
        three-dimensional. It follows that the center $\mathfrak{z}(\mathfrak{g})$
        of $\mathfrak{g}$ is at most three-dimensional, too. 
        \item If $\text{rank}\:\mathfrak{h}=2$, its Cartan subalgebra has to act as
        $\mathfrak{t}$ on $\text{Im}(\mathbb{O})$. The maximal trivial 
        $\mathfrak{h}$-submodule therefore is at most one-di\-men\-sional and we have 
        $\dim{\mathfrak{z}(\mathfrak{g})}\leq 1$. Moreover, $\text{rank}\:
        \mathfrak{g}$ has to be odd.
    \end{enumerate}

    \item \label{U1K} Let $G=G'\times U(1)$ and $H=H'\times U(1)$. If the second
    factors of both groups coincide, $U(1)$ acts trivially on $G/H$. If the
    second factor of $H$ is transversely embedded into $G'\times U(1)$, $G/H$ is
    covered by $G'/H'$. Since the group which acts on $G'/H'$ is $G'\times U(1)$
    instead of $G'$, we consider this case as a new one. The only other case
    which we have to consider is where $H'\times U(1) \subseteq G'$.

    \item Let $\mathfrak{m}$ be the orthogonal complement of $\mathfrak{h}$ in
    $\mathfrak{g}$ with respect to a biinvariant metric on $\mathfrak{g}$. The
    restriction of the adjoint action to a map $\mathfrak{h} \rightarrow
    \mathfrak{gl}(\mathfrak{m})$ is equivalent to the action of $\mathfrak{h}$
    on the tangent space of $G/H$. This identification helps us to compute the
    action of $\mathfrak{h}$ explicitly. In general, we omit that computation
    and give the reader a description of the isotropy action instead.
\end{enumerate}

In this section, we assume that $G=G'\times U(1)$ and $G/H = G'/H
\times S^1$. $G'/H$ admits a $G'$-invariant $SU(3)$-structure. We
can prove by similar arguments as in Lemma \ref{G2Hom1} and
\ref{G2Hom2} that our task reduces to classifying all
six-dimensional $G'$-homogeneous spaces $G'/H$ with $H \subseteq
SU(3)$. The possibilities for $\mathfrak{h}$ are thus fewer than
in the general situation. We prove our classification result, by 
considering each possible $\mathfrak{h}\subseteq \mathfrak{su}(3)$ 
separately. For reasons of brevity, we mostly mention only those 
$\mathfrak{g}$ which cannot be excluded by the above techniques. 

\underline{$\mathfrak{h}=\{0\}$:} In this case, $G/H$ simply is a
seven-dimensional, compact, connected Lie group. Up to coverings,
the only groups of this kind are $U(1)^7$, $SU(2)\times U(1)^4$,
and $SU(2)^2\times U(1)$.

\underline{$\mathfrak{h}=\mathfrak{u}(1)$:} Since
$\dim{\mathfrak{g}}=8$ and spaces of type $SU(3)/U(1)$ are
irreducible, the only remaining possibilities for $G$ are $SU(2)\times
U(1)^5$ and $SU(2)^2\times U(1)^2$. The first case can be
excluded, since the center of $G$ is too large. If
$G=SU(2)^2\times U(1)^2$, $H$ is embedded into $G$ by a map of
type:

\begin{equation}
e^{i\varphi} \mapsto \left(\,\left(\,\begin{array}{cc}
e^{ik_1\varphi} & 0 \\ 0 &
e^{-ik_1\varphi}\end{array}\,\right),\left(\,\begin{array}{cc}
e^{ik_2\varphi} & 0 \\ 0 & e^{-ik_2\varphi}\end{array}\,\right),
e^{ik_3\varphi},e^{ik_4\varphi}\,\right)\:,
\end{equation}

where $k_1,\ldots,k_4 \in\mathbb{Z}$. We repeat the argument
from page \pageref{U1K} twice and see that $G/H$ is
covered by $S^3\times S^3 \times S^1$ or that $H\subseteq
SU(2)^2$. The action of $H$ on the tangent space has at most two
non zero weights. We compare the weights of that action with the 
weights with which the one-dimensional subgroups of $G_2$ act on 
$\text{Im}(\mathbb{O})$. After that, we see that we can assume 
$|k_1|=|k_2|=1$. Since we obtain the same space for different choices
of the signs of $k_1$ and $k_2$, we can even assume that $k_1=k_2=1$.
If $(k_3,k_4)= (1,0)$, $G/H$ is diffeomorphic to $S^3\times S^3\times 
S^1$, and if $(k_3,k_4)= (0,0)$, we obtain the only space $G/H$ which
is not covered by $S^3\times S^3\times S^1$.

\underline{$\mathfrak{h}=\mathfrak{su}(2)$:} In this situation,
$G$ has to be a ten-dimensional compact Lie group. On the one
hand, $\dim{\mathfrak{z}( \mathfrak{g})}$ has to be positive,
since $G/H$ is $S^1$-reducible. On the other hand, we have
$\dim{\mathfrak{z}(\mathfrak{g})}\leq 3$. The only remaining
possibilities for $G$ therefore are $SU(2)^3\times U(1)$ and
$SU(3)\times U(1)^2$.

In the first case, we can embed $H$ diagonally, i.e. by the map
$g\mapsto (g,g,g,1)$. The action of $H$ on the tangent space is
the same as of $\mathfrak{su}(2)_{2,2}$ on $\text{Im}(\mathbb{O})$
and $G/H$ is diffeomorphic to $S^3\times S^3\times S^1$. If we had
embedded $H$ differently, it would act as the identity on a
four-dimensional subspace, which is impossible.

\label{EmbSU2SU3} In the second case, there are two possible
embeddings of $H$ into $SU(3)$: The first embedding is induced by
he standard representation of $SO(3)$ on
$\mathbb{R}^3\subseteq\mathbb{C}^3$. The only elements of $SU(3)$
which commute with all of $SO(3)$ are the multiples of the
identity. Since those elements are a discrete set, the action of
$H$ splits the tangent space into a trivial and a five-dimensional
irreducible submodule. There is no connected subgroup of $G_2$
which acts in this way on $\text{Im}(\mathbb{O})$ and we thus can
exclude this case. The second embedding is given by the following
map from $SU(2)$ to $SU(3)$:

\begin{equation}
A \mapsto \left(\,\begin{array}{cc} \hhline{-~}
\multicolumn{1}{|c|}{A} & \\ \hhline{--} & \multicolumn{1}{|c|}{1}
\\ \hhline{~-}
\end{array}\,\right)\:.
\end{equation}

In this situation, $H$ acts as $\mathbbm{V}_1^{\mathbb{C}} \oplus
3\mathbbm{V}_0^{\mathbb{R}}$ on the tangent space. Since
$\mathfrak{su}(2)_1$ acts in the same way, we have to put the
space $SU(3)/SU(2)\times U(1)^2 = S^5 \times T^2$ on our list. There
are no further embeddings of $H$ into $SU(3)$. This can be seen by 
considering the splitting of $\mathbb{C}^3$ into 
$\mathfrak{su}(2)$-submodules which is induced by the embedding of
$\mathfrak{su}(2)$ into $\mathfrak{su}(3)$.  

\underline{$\mathfrak{h}=2\mathfrak{u}(1)$:} Since
$\text{rank}\:\mathfrak{h}=2$ and $G/H$ is $S^1$-reducible, we have
$\dim{\mathfrak{z}(\mathfrak{g})}=1$. The group $G$ has to be
nine-dimensional. Therefore, we can assume that $G=SU(3)
\times U(1)$. Since $\mathfrak{su}(3) \subseteq \mathfrak{g}_2$
and $\text{rank} \:\mathfrak{su}(3) =
\text{rank}\:\mathfrak{g}_2$, any Cartan subalgebra of
$\mathfrak{su}(3)$ acts on $\mathbb{C}^3$ in the same way as
$\mathfrak{t}$ on $\text{span}(j,k,\ldots,k\epsilon)$. We thus
have to put the space $G/H = SU(3)/U(1)^2\times U(1)$ on our list.

\underline{$\mathfrak{h}=\mathfrak{su}(2)\oplus \mathfrak{u}(1)$:}
For similar reasons as in the previous case, $\mathfrak{g}$ has to
be the direct sum of a semisimple Lie algebra and
$\mathfrak{u}(1)$. With help of the classification of the
semisimple Lie algebras, we see that
$\mathfrak{g}=\mathfrak{sp}(2) \oplus \mathfrak{u}(1)$. We
describe a possible $G/H$ in detail. $Sp(2)$ has a subgroup
of type $Sp(1)\times U(1)$ which is given by:

\begin{equation}
H=\left\{\left(\,\begin{array}{cc} h_1 & 0 \\ 0 & h_2 \\
\end{array}\,\right)\:\left|
\begin{array}{l} \\ \\ \end{array}\right.
\!\!\!\!\! h_1\in\mathbb{H}, h_2\in\mathbb{C},
|h_1|=|h_2|=1\right\}\:.
\end{equation}

\label{su2u1arg}
The Lie algebra of $H$ acts in the same way on the tangent space
of $G/H$, which is diffeomorphic to $\mathbb{CP}^3\times S^1$, as
$\mathfrak{su}(2)_1 \oplus \mathfrak{u}(1)$ on $\text{Im}(\mathbb{O})$.
The kernel of the isotropy representation of $H$ is isomorphic to 
$\mathbb{Z}_2$. Therefore, we have an effective action by 
$(Sp(1)\times U(1))/\mathbb{Z}_2$ on the tangent space. Since that
group is isomorphic to $U(2)$, our example does not contradict the
fact that $G_2$ contains no subgroup of type $Sp(1)\times U(1)$.

We exclude the existence of further spaces of the above
kind. If $\mathfrak{g} = \mathfrak{sp}(2) \oplus \mathfrak{u}(1)$
and $\mathfrak{h} = \mathfrak{sp}(1)\oplus \mathfrak{u}(1)$,
either $G/H$ is covered by the sphere $Sp(2)/Sp(1)$, which is not
reducible, or $\mathfrak{h} \subseteq \mathfrak{sp}(2)$. 
\label{3su2salt} There are three embeddings of $\mathfrak{sp}(1)$
into $\mathfrak{sp}(2)$, which is isomorphic to $\mathfrak{so}(5)$.
In the first case, $\mathfrak{sp}(1)$ acts as $\mathfrak{so}(3)$ on
$\mathbb{R}^3 \subseteq \mathbb{R}^5$, in the second case, it acts as
$\mathfrak{su}(2)$ on $\mathbb{C}^2 \cong \mathbb{R}^4 \subseteq
\mathbb{R}^5$, and in the last case, it acts irreducibly on
$\mathbb{R}^5$. The second embedding yields the homogeneous space
$\mathbb{CP}^3\times S^1$, which we have described above. If
the semisimple part of $\mathfrak{h}$ was embedded by the first map,
it would act as $\mathfrak{su}(2)_{2,2}$ on the tangent space. Since
$\mathfrak{g}_2$ has no subalgebra of type $\mathfrak{su}(2)_{2,2}
\oplus \mathfrak{u}(1)$, this is not possible. It follows from Schur's
Lemma that there is no non zero element of $\mathfrak{so}(5)$ which 
commutes with the third possible embedding of the semisimple part. 
This case can therefore be excluded, too. 

\underline{$\mathfrak{h}=\mathfrak{su}(3)$:} As in the previous
cases, $G$ has to be a product of a 14-dimensional semisimple Lie
group $G'$ and $U(1)$. With help of the classification of the
semisimple Lie algebras, we conclude that $G'$ is $SU(3)\times
SU(2)^2$ or $G_2$. In the first case, $SU(3)$ acts trivially on
$G/H$ and in the second case we obtain $G_2/SU(3)\times U(1)$,
which is diffeomorphic to $S^6\times S^1$. We can verify that $H$
acts in the same way as the subgroup $SU(3)$ of $G_2$ on
$\text{Im}(\mathbb{O})$. Therefore, we have to put this space on
our list and have finally proven the first part of Theorem \ref{MainThm}.

\begin{Rem}
There is a one-to-one correspondence between the spaces from Theorem
\ref{MainThm}.\ref{RedTh} and the six-dimensional spaces which admit a
homogeneous $SU(3)$-structure. These spaces are considered by Cleyton 
and Swann \cite{Cley}, too. They obtain a list of homogeneous spaces 
which coincides with our list with the single exception of 
$SU(2)^2/U(1)\times T^2$, which seems to be missing in \cite{Cley}.
\end{Rem}

\section{The irreducible case}
\label{IrredSec}

In this section, we classify the $S^1$-irreducible spaces which admit a
homogeneous $G_2$-structure. As in the previous section, we
consider each possible $\mathfrak{h}$ separately.

\underline{$\mathfrak{h}=\{0\}$:} Since any seven-dimensional
compact Lie group is covered by a product of a semisimple Lie
group and a torus of positive dimension, we can exclude this case.

\underline{$\mathfrak{h}=\mathfrak{u}(1)$:} In the previous
section, we have already proven that if
$\mathfrak{h}=\mathfrak{u}(1)$ and $G/H$ is $S^1$-irreducible, we
necessarily have $G=SU(3)$. $G/H$ therefore is an Aloff-Wallach
space, i.e. a quotient $N^{k,l}:=SU(3)/U(1)_{k,l}$ with
$k,l\in\mathbb{Z}$ and

\begin{equation}
U(1)_{k,l}:=\left\{\left(\,\begin{array}{ccc} e^{ikt} & 0 & 0 \\ 0
& e^{ilt} & 0 \\ 0 & 0 & e^{-i(k+l)t}\\
\end{array}\,\right)\:\left|
\begin{array}{l} \\ \\ \\ \end{array}\right.
\!\!\!\!\! t\in\mathbb{R}\right\}\:.
\end{equation}

By an explicit calculation, we see that there exists a
one-dimensional Lie subalgebra of $\mathfrak{t}$ which acts in the
same way on $\text{Im}(\mathbb{O})$ as the Lie algebra of
$U(1)_{k,l}$ on the tangent space of $N^{k,l}$.

\underline{$\mathfrak{h}=\mathfrak{su}(2)$:} Since $\mathfrak{h}$
has to be embedded into the semisimple part of $\mathfrak{g}$,
$\mathfrak{z}(\mathfrak{g} )$ has to be trivial. Otherwise, $G/H$
would not be $S^1$-irreducible. The only remaining possibility for
$\mathfrak{g}$ therefore is $\mathfrak{so}(5)$. As we have
mentioned before, there are three embeddings of $\mathfrak{su}(2)$
into $\mathfrak{so}(5)$, which are distinguished by the splitting
of $\mathbb{R}^5$ with respect to $\mathfrak{su}(2)$:

\begin{enumerate}
    \label{3su2s}
    \item \underline{$\mathbb{R}^5=\mathbbm{V}_2^{\mathbb{R}} \oplus
    2\mathbbm{V}_0^{\mathbb{R}}$:} In this situation, $G/H$ is the
    Stiefel manifold $V^{5,2}=SO(5)/SO(3)$ of all orthonormal pairs in
    $\mathbb{R}^5$. The action of $\mathfrak{su}(2)$ splits the tangent
    space into $2\mathbbm{V}_2^{\mathbb{R}} \oplus
    \mathbbm{V}_0^{\mathbb{R}}$. Since $\mathfrak{su}(2)_{2,2}$ splits
    $\text{Im}(\mathbb{O})$ in the same way, $V^{5,2}$ admits an
    $SO(5)$-invariant $G_2$-structure.

    \item \underline{$\mathbb{R}^5=\mathbbm{V}_1^{\mathbb{C}} \oplus
    \mathbbm{V}_0^{\mathbb{R}}$:} If this is the case, $G/H$ is covered by
    the seven-sphere $Sp(2)/Sp(1)$. The action of $Sp(1)$ splits the tangent
    space into $\mathbbm{V}_1^{\mathbb{C}} \oplus 3\mathbbm{V}_0^{\mathbbm{R}}$.
    $\mathfrak{su}(2)_1$ acts in the same way on $\text{Im}(\mathbb{O})$ and
    $S^7$ thus admits an $Sp(2)$-invariant $G_2$-structure.

    \item \underline{$\mathbb{R}^5=\mathbbm{V}_4^{\mathbb{R}}$:} If
    $\mathfrak{su}(2)$ acts irreducibly on $\mathbb{R}^5$, it also acts
    irreducibly on the tangent space of $G/H$. Since the action of
    $\mathfrak{su}(2)_6$ on $\text{Im}(\mathbb{O})$ is irreducible, too, we have
    found another space which admits a homogeneous $G_2$-structure, namely the
    seven-dimensional Berger space $B^7$.
\end{enumerate}

\underline{$\mathfrak{h}=2\mathfrak{u}(1)$:} Since $\mathfrak{h}$
is of rank $2$, $\dim{\mathfrak{z}(\mathfrak{g})}$ is either $0$ or
$1$. If the center is one-dimensional, we have
$\mathfrak{g}=\mathfrak{su}(3)\oplus \mathfrak{u}(1)$ and
$\mathfrak{h}$ is transversely embedded into that direct sum. In
this situation, $G/H$ is covered by an Aloff-Wallach space
$N^{k,l}$, on which $SU(3) \times U(1)$ acts transitively. The
group $SU(3)$ acts as usual by left multiplication on $N_{k,l}$.
Moreover, a certain one-dimensional subgroup of the normalizer
$\text{Norm}_{SU(3)}U(1)_{k,l}$ acts on $N^{k,l}$ by
right multiplication. This subgroup can be identified with the
second factor of $SU(3)\times U(1)$.

If $\mathfrak{g}$ is semisimple, we can assume that $\mathfrak{g}
= 3\mathfrak{su}(2)$. We describe the possible embeddings of
$2\mathfrak{u}(1)$ into $3\mathfrak{su}(2)$. A Cartan subalgebra
of $3\mathfrak{su}(2)$ is given by:

\begin{equation}
\left\{\left(\,\begin{array}{cccccc} \hhline{--~~~~}
\multicolumn{1}{|c}{ix} & \multicolumn{1}{c|}{0} &&&& \\
\multicolumn{1}{|c}{0} & \multicolumn{1}{c|}{-ix} &&&& \\
\hhline{----~~} && \multicolumn{1}{|c}{iy} &
\multicolumn{1}{c|}{0} && \\ && \multicolumn{1}{|c}{0} &
\multicolumn{1}{c|}{-iy} && \\ \hhline{~~----} &&&&
\multicolumn{1}{|c}{iz} & \multicolumn{1}{c|}{0} \\ &&&&
\multicolumn{1}{|c}{0} & \multicolumn{1}{c|}{-iz} \\
\hhline{~~~~--}
\end{array}\,\right)\:\left|
\begin{array}{l} \\ \\ \\ \\ \\ \\ \end{array}
\right.\!\!\!\!\! x,y,z\in\mathbb{R}\right\}
\end{equation}

We fix the biinvariant metric $q(X,Y):=-\text{tr}(XY)$ on
$3\mathfrak{su}(2)$. Let $\mathfrak{k}_{k,l,m}$, where
$k,l,m\in\mathbb{Z}$, be the one-dimensional subalgebra of
$3\mathfrak{su}(2)$ which is generated by the matrix with $x=k$,
$y=l$, and $z=m$. Furthermore, let $2\mathfrak{u}(1)_{k,l,m}$ be
the $q$-orthogonal complement of $\mathfrak{k}_{k,l,m}$ in the above
Cartan subalgebra. Any connected two-dimensional Lie subgroup of
$SU(2)^3$ is conjugate to a connected subgroup with a Lie algebra
of type $2\mathfrak{u}(1)_{k,l,m}$. We denote the quotient of
$SU(2)^3$ by that subgroup by $Q^{k,l,m}$.

By the action of the group $(\mathbb{Z}_2)^3\rtimes S_3$ of outer
automorphisms of $3\mathfrak{su}(2)$, we can change the signs and
the order of $(k,l,m)$ arbitrarily. We may therefore assume
without loss of generality that $k\geq l \geq m\geq 0$. The
isotropy representation of $2\mathfrak{u}(1)_{k,l,m}$ on the
tangent space of $Q^{k,l,m}$ is with respect to a suitable basis
given by:

\begin{equation}
\label{QklmIsotropy}
\left\{\left(\,\begin{array}{ccccccc} \hhline{--~~~~~}
\multicolumn{1}{|c}{0} & \multicolumn{1}{c|}{x} &&&&& \\
\multicolumn{1}{|c}{-x} & \multicolumn{1}{c|}{0} &&&&& \\
\hhline{----~~~} && \multicolumn{1}{|c}{0} &
\multicolumn{1}{c|}{y} &&& \\ && \multicolumn{1}{|c}{-y} &
\multicolumn{1}{c|}{0} &&& \\ \hhline{~~----~} &&&&
\multicolumn{1}{|c}{0} & \multicolumn{1}{c|}{z} & \\ &&&&
\multicolumn{1}{|c}{-z} & \multicolumn{1}{c|}{0} & \\
\hhline{~~~~---} &&&&&& \multicolumn{1}{|c|}{0} \\
\hhline{~~~~~~-}
\end{array}\,\right)\:\left|
\begin{array}{l} \\ \\ \\ \\ \\ \\ \\ \end{array}
\right.\!\!\!\!\! xk+yl+zm=0\right\}
\end{equation}

By comparing (\ref{QklmIsotropy}) with the Cartan subalgebra
(\ref{csag2}) of $\mathfrak{g}_2$, we see that only $Q^{1,1,1}$ 
admits an $SU(2)^3$-invariant $G_2$-structure.

\underline{$\mathfrak{h}=\mathfrak{su}(2)\oplus \mathfrak{u}(1)$:}
Since $\text{rank} \: \mathfrak{su}(2)\oplus \mathfrak{u}(1)=2$,
the center of $\mathfrak{g}$ is at most one-dimensional.
$\mathfrak{g}$ has to be an eleven-dimensional Lie algebra and
therefore is either $\mathfrak{su}(3) \oplus \mathfrak{su}(2)$ or
$\mathfrak{so}(5)\oplus\mathfrak{u}(1)$.

We start with the first of the two cases. The semisimple part of
$\mathfrak{h}$ we denote by $\mathfrak{su}(2)'$. In order to
classify the homogeneous spaces which we can obtain in this
situation, we have to describe the possible embeddings of 
$\mathfrak{su}(2)'$ into $\mathfrak{su}(3)\oplus\mathfrak{su}(2)$.
$\mathfrak{su}(2)'\cap\mathfrak{su}(3)$ has to be nontrivial.
Otherwise, $G$ would not act almost effectively on $G/H$. The
projection of $\mathfrak{su}(2)'$ onto $\mathfrak{su}(3)$
therefore has to be one of the two maps which we have described on
page \pageref{EmbSU2SU3}. If $\mathfrak{su}(2)'$ acted irreducibly
on $\mathbb{C}^3$, the tangent space of $G/H$ would contain a
five-dimensional $\mathfrak{su}(2)'$-submodule. This follows by
the same arguments as on page \pageref{EmbSU2SU3}. Since no
subalgebra of $\mathfrak{g}_2$ acts in this way on
$\text{Im}(\mathbb{O})$, $\mathfrak{su}(2)'$ has to split $\mathbb{C}^3$
into $\mathbbm{V}_1^{\mathbb{C}} \oplus \mathbbm{V}_0^{\mathbb{C}}$.
Next, we consider the projection of $\mathfrak{su}(2)'$ onto the
second summand of $\mathfrak{su}(3) \oplus \mathfrak{su}(2)$. We
first assume that $\mathfrak{su}(2)' \subseteq \mathfrak{su}(3)$.
In this situation, the center of $\mathfrak{h}$ is without loss of
generality generated by a matrix of type

\begin{equation}
\left(\,\begin{array}{ccccc} \hhline{---~~}
\multicolumn{1}{|c}{ki} & 0 & \multicolumn{1}{c|}{0} && \\
\multicolumn{1}{|c}{0} & ki & \multicolumn{1}{c|}{0} && \\
\multicolumn{1}{|c}{0} & 0 & \multicolumn{1}{c|}{-2ki} && \\
\hhline{-----} &&& \multicolumn{1}{|c}{li} &
\multicolumn{1}{c|}{0} \\ &&& \multicolumn{1}{|c}{0} &
\multicolumn{1}{c|}{-li} \\ \hhline{~~~--}
\end{array}\,\right)
\end{equation}

where $k$ and $l$ are integers. $\mathfrak{su}(2)'$ acts as
$\mathfrak{su}(2)_1$ on the tangent space of $G/H$. There is up to
conjugation only one one-dimensional subalgebra of
$\mathfrak{g}_2$ which commutes with $\mathfrak{su}(2)_1$.
Therefore, the weights with which the center of $\mathfrak{h}$
acts on the tangent space are uniquely determined. By computing
the action of the above matrix on the tangent space, we see that
we necessarily have $l=\pm 3k$. The quotient $G/H$ is in both
cases up to an $SU(3)\times SU(2)$-equivariant diffeomorphism the
same and admits an $SU(3)\times SU(2)$-invariant $G_2$-structure.
We use the same notation as Castellani \cite{Cast} and call our
space $M^{1,1,0}$.

If the projection of $\mathfrak{su}(2)'$ onto the second summand
of $\mathfrak{su}(3) \oplus \mathfrak{su}(2)$ is bijective, there
is up to conjugation only one one-dimensional subalgebra of
$\mathfrak{su}(3) \oplus \mathfrak{su}(2)$ which commutes with
$\mathfrak{su}(2)'$. In this situation, $G/H$ is diffeomorphic to
the exceptional Aloff-Wallach space $N^{1,1}$. $SU(3)$ acts on a
$g U(1)_{1,1}$ by matrix multiplication from the left. Since
$U(1)_{1,1}$ commutes with $S(U(2)\times U(1))$ which is
isomorphic to $SU(2)$, $g U(1)_{1,1}\mapsto gh^{-1} U(1)_{1,1}$
defines a left action by $SU(2)$ on $N^{1,1}$ which commutes with
the action of $SU(3)$. The isotropy group of the $SU(3)\times
SU(2)$-action which we have defined is $SU(2)\times U(1)$. The
embedding of its Lie algebra into $\mathfrak{su}(3) \oplus
\mathfrak{su}(2)$ is the same as we have described above. We thus
have found another group action on $N^{1,1}$ which we have to
include in our list.

In both of the above two cases, there exists an element of $G$
which is of order two and acts trivially on $G/H$. For the same
reasons as on page \pageref{su2u1arg}, the fact that $G_2$ contains
no subgroup of type $SU(2) \times U(1)$ therefore does not contradict
the statement of our theorem.

Next, we assume that $\mathfrak{g}=\mathfrak{so}(5) \oplus
\mathfrak{u}(1)$. The embedding of $\mathfrak{su}(2)'$ into
$\mathfrak{so}(5)$ has to be one of the three subalgebras which we
have described on pages \pageref{3su2salt} and \pageref{3su2s}.
Furthermore, the projection of $\mathfrak{z}( \mathfrak{h})$ onto
$\mathfrak{so}(5)$ should not be trivial. If $\mathfrak{su}(2)'$
was embedded by its five-dimensional representation into
$\mathfrak{so}(5)$, there would be no element of
$\mathfrak{so}(5)$ left which commutes with $\mathfrak{su}(2)'$.
Since this contradicts our statement on
$\mathfrak{z}(\mathfrak{h})$, we can exclude this case. If
$\mathfrak{su}(2)'$ was embedded by its three-dimensional
representation, it would decompose its complement in
$\mathfrak{so}(5)$ into $2\mathbbm{V}_2^{\mathbb{R}} \oplus
\mathbbm{V}_0^{\mathbbm{R}}$. $\mathfrak{g}_2$ has no subalgebra of
type $\mathfrak{su}(2)_{2,2} \oplus \mathfrak{u}(1)$ and we thus can
exclude this case, too. The only remaining case is where 
$\mathfrak{su}(2)'$ is embedded by its two-dimensional complex 
representation. Since $\mathfrak{z}(\mathfrak{h})$ has to commute with
$\mathfrak{su}(2)'$, its projection onto $\mathfrak{so}(5)$ has to
be an element of the second summand of $\mathfrak{su}(2)' \oplus
\mathfrak{su}(2)$, which we identify with the Lie subalgebra
$\mathfrak{so}(4)$ of $\mathfrak{so}(5)$. If
$\mathfrak{h}\subseteq \mathfrak{so}(5)$, we obtain the space
$\mathbb{CP}^3\times S^1$, which we already have considered in the
previous section. If this is not the case, $G/H$ is covered by
$S^7$, which is equipped with an action of $Sp(2)\times U(1)$.

\underline{$\mathfrak{h}=2\mathfrak{su}(2)$:} Since
$\dim{\mathfrak{h}}=6$, the dimension of $\mathfrak{g}$ has to be
$13$. There is no non zero element of $\text{Im}(\mathbb{O})$ on
which the subalgebra $2\mathfrak{su}(2)$ of $\mathfrak{g}_2$ acts
trivially. Therefore, $\mathfrak{z}(\mathfrak{g})$ has to be
trivial. The only remaining possibility for $\mathfrak{g}$ is
$\mathfrak{so}(5) \oplus \mathfrak{su}(2)$.

It follows from Lemma \ref{G2Hom1} and Theorem \ref{G2SubThm} that
$\mathfrak{h}$ has to decompose the tangent space into
$\mathbbm{V}_{2,0}^{\mathbb{R}} \oplus
\mathbbm{V}_{1,1}^{\mathbb{C}}$. Let $\imath:2\mathfrak{su}(2)
\rightarrow \mathfrak{so}(5) \oplus \mathfrak{su}(2)$ be the
embedding of $\mathfrak{h}$ into $\mathfrak{g}$, $\pi_1:
\mathfrak{so}(5) \oplus \mathfrak{su}(2) \rightarrow
\mathfrak{so}(5)$ be the projection on the first summand, and
$\pi_2:\mathfrak{so}(5) \oplus \mathfrak{su}(2) \rightarrow
\mathfrak{su}(2)$ be the projection on the second one. The tangent 
space contains a submodule of type $\mathbbm{V}_{1,1}^{\mathbb{C}}$
only if $(\pi_1\circ\imath)(2\mathfrak{su}(2))$ is the standard
embedding of $\mathfrak{so}(4)$ into $\mathfrak{so}(5)$. The first
summand of $2\mathfrak{su}(2)$ has to act irreducibly on a
three-dimensional submodule of the tangent space and we therefore
can assume that

\begin{equation}
(\pi_2\circ\imath)(x,y)=x \quad \forall x,y\in\mathfrak{su}(2)\:.
\end{equation}

We are now able to describe $G/H$ explicitly. Let $S^7\subseteq
\mathbb{H}^2$ be the seven-sphere. $Sp(2)$ acts on $S^7$ from the
left by matrix multiplication. We identify $Sp(1)$ with the group
of all unit quaternions. Since the scalar multiplication on a
quaternionic vector space acts from the right,
scalar multiplication with $h^{-1}$ where $h\in Sp(1)$ defines a
left action of $Sp(1)$ on $S^7$. We thus have constructed a
transitive $Sp(2)\times Sp(1)$-action on $S^7$. The isotropy group
of this action is $Sp(1)\times Sp(1)$ and the isotropy action has
the properties which we have demanded above. Analogously to the
case where $H=SU(2)\times U(1)$, the kernel of the isotropy representation
of $Sp(1)\times Sp(1)$ is $\mathbb{Z}_2$ and the group which acts
effectively on the tangent space is in fact $(Sp(1)\times Sp(1))/
\mathbb{Z}_2$, which is isomorphic to $SO(4)$.

\underline{$\mathfrak{h}=\mathfrak{su}(3)$:} $G$ has to be a Lie
group of dimension $15$ which contains $SU(3)$. With help of the
classification of the compact Lie groups, we see that $G$ is
covered either by a product of $SU(3)$ and a seven-dimensional Lie
group or by $SU(4)$. In the first case, $G$ would not act almost
effectively on $G/H$. In the second case, $G/H$ is covered by
$S^7$.

\underline{$\mathfrak{h}=\mathfrak{g}_2$:} For similar reasons as
above, we have $\mathfrak{g}=\mathfrak{so}(7)$. Therefore, $G/H$
is covered by the seven-dimensional sphere $\text{Spin}(7)/G_2$ and
we have completed the proof of Theorem \ref{MainThm}.\ref{IrredTh}.

\begin{Rem}
Friedrich, Kath, Moroianu, and Semmelmann \cite{Fried} have classified
all spaces which admit a homogeneous nearly parallel $G_2$-structure. In
particular, the authors prove that all spaces from Theorem 
\ref{MainThm}.\ref{IrredTh} admit such a $G_2$-structure. In the table of our
theorem, we also have listed all transitive group actions on those spaces
which preserve a $G_2$-structure. On the sphere $S^7$, for example, there are 
$G_2$-structures which are invariant under $\text{Spin}(7)$, $SU(4)$,
$Sp(2)\times Sp(1)$, $Sp(2)\times U(1)$, or $Sp(2)$. We remark that
some of the Aloff-Wallach spaces are diffeomorphic or homeomorphic 
to each other, although they are not $SU(3)$-equivariantly diffeomorphic.
This phenomenon is discussed by Kreck and Stolz \cite{Kreck}.
Their results prove that on the same space there can exist
$G_2$-structures which are preserved by different transitive Lie
group actions.
\end{Rem}

\section{Existence of the cosymplectic $G_2$-structures}
\label{CosympSec}

In the previous two sections, we have classified all spaces which
admit a homogeneous $G_2$-structure. The aim of this section is to
prove that a transitive group action which leaves at least one
$G_2$-structure invariant also leaves a cosymplectic
$G_2$-structure invariant. We prove this fact by a case-by-case
analysis. Although most of this work has already been done by
other $\text{authors}$, there are still some cases left open.

Since any nearly parallel $G_2$-structure is also cosymplectic,
the article of Friedrich et al. \cite{Fried} answers our question
for many subcases of the irreducible case. More precisely, we only
have to consider those irreducible spaces on which we have more
than one transitive group action.

Let $S^7\subseteq \mathbb{O}$ be the unit sphere. We equip the
tangent space $\text{Im}(\mathbb{O})$ of $1\in\mathbb{O}$ with the
canonical $G_2$-structure $\omega$ from page \pageref{omega}. By
the action of $\text{Spin}(7)$, we can extend $\omega$ to a nearly
parallel $G_2$-structure on all of $S^7$. Since we have
$Sp(2)\subseteq SU(4)\subseteq \text{Spin}(7)$, $\omega$ is
invariant with respect to the action of the three groups. In
\cite{Fried}, the authors describe a homogeneous nearly parallel
$G_2$-structure on $S^7$. The associated metric on $S^7$ is the 
squashed one and its isometry group is $Sp(2)\times Sp(1)$. Since
the $G_2$-structure is homogeneous, it has to be at least 
$Sp(2)$-invariant. We assume that the second factor of $Sp(2) 
\times Sp(1)$ does not preserve the $G_2$-structure. In that 
situation, there exists a one-dimensional subgroup of $Sp(1)$ which
generates a continuous family of nearly parallel $G_2$-structures but
preserves the associated metric. Any nearly parallel $G_2$-structure 
induces a Killing spinor and the dimension of the space of all Killing
spinors thus is at least two. Since it is known (cf. \cite{Fried}) 
that this dimension is in fact one, the $G_2$-structure is $Sp(2)\times
Sp(1)$- and in particular $Sp(2)\times U(1)$-invariant. All in all, we
have found for each transitive action on $S^7$ an invariant cosymplectic
$G_2$-structure.

Next, we consider the Aloff-Wallach spaces. Cveti\v{c} et al.
\cite{Cve} have proven that any Aloff-Wallach space admits two
$SU(3)$-invariant nearly parallel $G_2$-structures, which coincide
for $k=-l$. It is known (cf. \cite{Fried}) that the isometry group
of the associated metric is $SU(3)\times U(1)$. Since the space of all
Killing spinors is one-dimensional (cf. \cite{Cve}, \cite{Fried}), we
can conclude by the same arguments as above that both $G_2$-structures 
are not only $SU(3)$- but also $SU(3)\times U(1)$-invariant.

The nearly parallel $G_2$-structure on $N^{1,1}$ which is
considered in \cite{Fried} is preserved by $SU(3)\times SU(2)$.
Since $SU(3) \subseteq SU(3)\times U(1) \subseteq SU(3)\times
SU(2)$, that $G_2$-structure is invariant with respect to all of
the three group actions from Theorem \ref{MainThm}.\ref{IrredTh}.

We proceed to the reducible case. Butruille \cite{Butruille} has
proven that the only six-dimensional manifolds which admit a
homogeneous nearly K\"ahler structure are $S^6$, $\mathbb{CP}^3$,
$SU(3)/U(1)^2$, and $S^3\times S^3$. These four manifolds have also
been considered by B\"ar \cite{Baer}, since they carry a real Killing 
spinor. The groups which preserve the nearly K\"ahler structure 
on the first three spaces are $G_2$, $Sp(2)$, and $SU(3)$. In 
\cite{Baer} it is also proven that $S^3\times S^3$ admits a 
nearly K\"ahler structure which is invariant under an 
$SU(2)^3$-action. The isotropy group of this action is $SU(2)$, 
which is embedded as the diagonal subgroup by

\begin{equation}
g\mapsto (g,g,g)\:.
\end{equation}

We denote the metric, the real two-form, and the complex
$(3,0)$-form which determine the $SU(3)$-structure by $g$,
$\alpha$, and $\theta$. Furthermore, we denote the real
(imaginary) part of $\theta$ by $\theta^{\text{Re}}$
($\theta^{\text{Im}}$). We have $d\alpha=
3\lambda\:\theta^{\text{Re}}$ and $d\theta^{\text{Im}}=
-2\lambda\: \alpha \wedge \alpha$ for a $\lambda\in\mathbb{R}
\setminus\{0\}$, since the four spaces are nearly K\"ahler. These
equations are discussed in more detail by Hitchin \cite{Hitchin}.
On a product of a circle and a nearly K\"ahler manifold of real
dimension six, we can define a $G_2$-structure by
$\omega:=\alpha\wedge dt + \theta^{Im}$. Here, "$t$" denotes the
coordinate of the circle. By a straightforward calculation, it
follows that $d \ast\omega=0$. All in all, we have proven our
statement for the last three spaces from Theorem 
\ref{MainThm}.\ref{RedTh} and for all three actions on $S^3 
\times S^3 \times S^1$.

On the torus $T^7$, we have the flat $G_2$-structure, which is of
course cosymplectic. On $\mathbb{C}^2\times T^4$ ($\mathbb{C}^3
\times T^2$), there exists a flat $\text{Spin}(7)$-structure
$\Omega$. It is preserved by the action of $SU(2)\times U(1)^4$
($SU(3) \times U(1)^2$), where the first factor acts on
$\mathbb{C}^2$ ($\mathbb{C}^3$) and the second one by 
translations on the torus. The principal orbits of this action,
which is of cohomogeneity one, are $S^3\times T^4$ ($S^5\times T^2$).
$\Omega$ induces an $SU(2)\times U(1)^4$ ($SU(3) \times U(1)^2$)-invariant
$G_2$-structure on any principal orbit. This $G_2$-structure is
cosymplectic, since $d\Omega=0$.

The only remaining space is $SU(2)^2/U(1)\times T^2$. The issue of
homogeneous $G_2$-structures on this space is not yet discussed in
the literature. In the following, we construct an explicit
$SU(2)^2\times U(1)^2$-invariant cosymplectic $G_2$-structure on
$SU(2)^2/U(1)\times T^2$. First, we choose the following basis of
$\mathfrak{su}(2)$:

\begin{equation}
\sigma_1:=\left(\,\begin{array}{cc} i & 0 \\ 0 & -i \\
\end{array}\,\right)\:,\quad
\sigma_2:=\left(\,\begin{array}{cc} 0 & -1 \\ 1 & 0 \\
\end{array}\,\right)\:,\quad
\sigma_3:=\left(\,\begin{array}{cc} 0 & i \\ i & 0 \\
\end{array}\,\right)\:.
\end{equation}

The above basis obeys the following commutator relations:

\begin{equation}
[\sigma_1,\sigma_2]:=-2\sigma_3\:,\quad
[\sigma_2,\sigma_3]:=-2\sigma_1\:,\quad
[\sigma_3,\sigma_1]:=-2\sigma_2\:.
\end{equation}

As usual, the tangent space of $SU(2)^2/U(1)\times T^2$ can be
identified with the complement $\mathfrak{m}$ of the isotropy
algebra in $2\mathfrak{su}(2) \oplus 2\mathfrak{u}(1)$. We
construct a basis $(e_1,\ldots,e_7)$ of $\mathfrak{m}$ and
supplement it with a generator $e_8$ of $\mathfrak{h}$ to a basis
of $2\mathfrak{su}(2)\oplus 2\mathfrak{u}(1)$. Let $e_1$ and $e_2$
be generators of the center of $2\mathfrak{su}(2) \oplus 
2\mathfrak{u}(1)$. Furthermore, let

\begin{equation}
\begin{array}{lll}
e_3:=\left(\,\begin{array}{cc} \hhline{-~}
\multicolumn{1}{|c|}{\sigma_1} & \\ \hhline{--} &
\multicolumn{1}{|c|}{-\sigma_1} \\ \hhline{~-}
\end{array}\,\right) \:, &
e_4:=\left(\,\begin{array}{cc} \hhline{-~}
\multicolumn{1}{|c|}{\sigma_2} & \\ \hhline{--} &
\multicolumn{1}{|c|}{0} \\ \hhline{~-}
\end{array}\,\right) \:, &
e_5:=\left(\,\begin{array}{cc} \hhline{-~}
\multicolumn{1}{|c|}{\sigma_3} & \\ \hhline{--} &
\multicolumn{1}{|c|}{0} \\ \hhline{~-}
\end{array}\,\right)\:, \\
&& \\ e_6:=\left(\,\begin{array}{cc} \hhline{-~}
\multicolumn{1}{|c|}{0} & \\ \hhline{--} &
\multicolumn{1}{|c|}{\sigma_2} \\ \hhline{~-}
\end{array}\,\right) \:, &
e_7:=\left(\,\begin{array}{cc} \hhline{-~} \multicolumn{1}{|c|}{0}
& \\ \hhline{--} & \multicolumn{1}{|c|}{\sigma_3} \\ \hhline{~-}
\end{array}\,\right) \:, &
e_8:=\left(\,\begin{array}{cc} \hhline{-~}
\multicolumn{1}{|c|}{\sigma_1} & \\ \hhline{--} &
\multicolumn{1}{|c|}{\sigma_1} \\ \hhline{~-}
\end{array}\,\right)\:. \\
\end{array}
\end{equation}

We define by $e^i(e_j):=\delta^i_j$ a basis of left invariant
one-forms on $SU(2)^2\times U(1)^2$. With help of the formula
$de^i(e_j,e_k) = - e^i([e_j,e_k])$ and the commutator relations on
$2\mathfrak{su}(2)\oplus 2\mathfrak{u}(1)$ it follows that:

\begin{equation}
\label{MaurerCartan}
\begin{array}{rcl}
de^1 & = & 0 \\ de^2 & = & 0 \\ de^3 & = & e^4\wedge e^5 -
e^6\wedge e^7 \\ de^4 & = & -2 e^3\wedge e^5 + 2 e^5\wedge e^8\\
de^5 & = & 2 e^3\wedge e^4 - 2 e^4\wedge e^8\\ de^6 & = & 2
e^3\wedge e^7 + 2 e^7\wedge e^8\\ de^7 & = & -2 e^3\wedge e^6 - 2
e^6\wedge e^8\\ de^8 & = & e^4 \wedge e^5 + e^6 \wedge e^7 \\
\end{array}
\end{equation}

In order to construct a homogeneous cosymplectic $G_2$-structure, we
introduce a further basis $(f_1,\ldots,f_7)$ of $\mathfrak{m}$: 

\begin{equation}
\begin{array}{lll}
f_1:=e_1\:, & f_2:=e_2\:, & f_3:=e_3\:, \\
f_4:= e_5 + e_7\:, & f_5:= -e_4-e_6\:,\\
f_6:= e_5 - e_7\:, & f_7:= -e_4 + e_6\:. \\ 
\end{array}
\end{equation}

With respect to this basis, the action of $e_8$ on $\mathfrak{m}$
is represented by a matrix which is contained in the Cartan 
subalgebra (\ref{csag2}). Therefore, we can identify
$(f_1,\ldots,f_7)$ with the standard basis of
$\text{Im}(\mathbb{O})$. This identification yields an $SU(2)^2\times
U(1)^2$-invariant $G_2$-structure $\omega$ on $SU(2)^2/U(1)\times
T^2$, which satisfies:

\begin{equation}
\ast\omega=-2e^{1245} + 2e^{1267} -2e^{1346} -2e^{1357} - 2e^{2347} + 2e^{2356} + 
4e^{4567}\:.
\end{equation}

As in Section \ref{2ndSec}, $e^{ijkl}$ is an abbreviation of 
$e^i \wedge e^j \wedge e^k \wedge e^l$. With help of the equations 
(\ref{MaurerCartan}) and the fact that the projection of $e_8$ onto
$SU(2)^2/U(1)\times T^2$ vanishes, we are able to compute 
$d\ast\omega$  and see that our $G_2$-structure is indeed cosymplectic. This 
calculation finishes the proof of Theorem \ref{MainThm}.

\begin{Rem}
Our proof that any $G$-homogeneous space $G/H$ which admits an arbitrary
$G$-invariant $G_2$-structure also admits a cosymplectic one is
done by a case-by-case analysis. If $G/H$ is irreducible, there
even exists a $G$-invariant nearly parallel $G_2$-structure on
$G/H$. The author suspects that it is possible to prove these
facts more directly.
\end{Rem}

\end{document}